\documentclass[11pt,twoside]{amsart}
\usepackage{latexsym}
\usepackage{amssymb,amsbsy,amsmath,amsfonts,amssymb,amscd}
\usepackage{graphicx}
\usepackage{hyperref}
\usepackage{xcolor}
\usepackage{dsfont}

\newcommand{\commentout}[1]{}

\newcommand{\R}{\mathbb{R}}
\newcommand{\N}{\mathbb{N}}

\newcommand {\Chi} {{\bf \raise 2pt \hbox{$\chi$}} }

\newcommand {\ep}  {\epsilon}
\newcommand{\beq}{\begin{equation}}
\newcommand{\eeq}{\end{equation}}
\newcommand{\bea} {\begin{array}{rl}}
\newcommand{\eea} {\end{array}}
\newcommand{\bepa}{\left\{ \begin{array}{l}}
\newcommand{\eepa} {\end{array}\right.}
\newtheorem{theorem}{Theorem}[section]
\newtheorem{lemma}[theorem]{Lemma}

\newtheorem{proposition}[theorem]{Proposition}

\def \Ž{\'e}
\def \ˆ{\`a}
\def \{\`e}
\def \{\^e}
\def \{\c{c}}
\def \"{\^{\i}}
\def \™{\^o}
\def \ž{\^u}
\def \‰{\^a}
\def \{\`u}
\def \Š{\"a}
\def\'{\"e}
\def\š{\"o}
\def\•{\"i}
\def\Ÿ{\"u}
\def\…{\"O}
\newcommand {\no}{\noindent}

\title{ {New regularity results and long time behavior \\
  of pathwise (stochastic) Hamilton-Jacobi equations}}

 \author{
Pierre-Louis Lions$^{1,3}$ and Panagiotis E. Souganidis$^{2,4}$}

\begin{document}

\maketitle
\pagestyle{plain}

\pagenumbering{arabic}

\begin{abstract}

%\smallskip
\no We present two new sharp regularity results (regularizing effect and propagation of regularity) for viscosity solutions of uniformly convex space homogeneous Hamilton-Jacobi equations. In turn, these estimates yield new intermittent   stochastic regularization results for pathwise  (stochastic) viscosity solutions of Hamilton-Jacobi equations with uniformly convex Hamiltonians  and rough multiplicative time dependence. Finally, we use the intermittent estimates to study the long time behavior of the pathwise (stochastic) viscosity solutions of convex Hamilton-Jacobi equations. 
\end{abstract}

\medskip

\\
\\
\noindent {\bf AMS Class. Numbers:} 60H15, 35D40.
\bigskip
\section{Introduction}
\subsection*{The general setting and results} 
We present two new sharp regularity results (regularizing effect and propagation of regularity) for viscosity solutions of 
\begin{equation}\label{takis1}
u_t=\pm H(Du) \ \ \text{in} \ \  \R^d\times (0,\infty), \end{equation}
where 
\begin{equation}\label{hconvex2}
H\in C^2(\R^d) \ \ \text{is uniformly convex}, %\ \ and} \ \ H(p)>H(0)=0 \ \ \text{for all} \ \ p\in \R^d\setminus\{0\},
\eeq
%and
%\begin{equation}\label{takis1.1}
%H(p)>H(0)=0 \ \ \text{for all} \ \ p\in \R^d\setminus\{0\}.
%\end{equation}
It follows from \eqref{hconvex2} that is, there  exist $\Theta, \theta >0$ such that, for all $p \in \R^d$ and in the sense of symmetric matrices,
\begin{equation}\label{hconvex1}
\theta I \leq D^2 H(p) \leq \Theta I,
%\Theta |\xi|^2 \geq \sum_{i,j=1}^d H_{p_i p_j} \xi_i\xi_j \geq \theta |\xi|^2.
\end{equation}
where $I$ is the identity matrix in $\R^d$.  The upper upper in \eqref{hconvex1} can be relaxed when dealing  with Lipshitz continuous solutions of \eqref{takis1}.  
\smallskip

%It follows from \eqref{hconvex2} that $H$ has a strict minimum. For simplicity we assume 

%Finally, \eqref{takis1.1} means that $H$ has a strict minimum. The location of this minim
%\smallskip]

The new estimates are then used to obtain two new stochastic regularization-type results (Lipshitz continuity and $C^{1,1}$-regularity, the latter only when $d=1$) for pathwise (stochastic) viscosity solutions with uniformly convex Hamiltonians and rough multiplicative time dependence, that is, 
\begin{equation}\label{takis100}
du=H(Du)\cdot d\zeta \ \ \text{in} \ \  \R^d\times (0,\infty),
% \quad \text{and} \quad  u(\cdot, 0)=u_0 \ \ \text{in } \ \ \R^d,
\end{equation}
with  $\zeta \in C_0([0,\infty);\R)$, the set of continuous functions $\zeta:[0,\infty)\to \R$ such that $\zeta(0)=0.$
\smallskip

Finally we use the new estimates to  investigate the long time behavior of periodic pathwise solutions of \eqref{takis100}, when %  for  Hamiltonians such that 
%
%
%that is, when 
%
% and rough multiplicative time dependence, that is, 
%\begin{equation}\label{takis100}
%du=H(Du)\cdot d\zeta \ \ \text{in} \ \  \R^d\times (0,\infty) \quad \text{and} \quad  u(\cdot, 0)=u_0 \ \ \text{in } \ \ \R^d,
%\end{equation}
%with  $\zeta \in C_0([0,\infty);\R)$, the set of continuous functions $\zeta:[0,\infty)\to \R$ such that $\zeta(0)=0,$  and 
\begin{equation}\label{takis712}
H\in C(\R^d) \ \text{is convex \ \ and } \ \ H(p)>H(0)=0 \ \ \text{for all} \ \ p\in \R^d\setminus\{0\}. % \ \ \text{convex \quad and} \quad H(p) > H(0)=0 \quad \text{for all } \quad p\in \R^d\setminus \{0\}.
\end{equation}

%and $u_0$ periodic in the unit cube. 
\smallskip

%We also present two new surprising and sharp regularizing results for the viscosity solutions of the ``deterministic'' Hamilton-Jacobi equation
%\begin{equation}\label{takis1.1}
%u_t=H(Du) \ \ \text{in} \ \  \R^d\times (0,\infty).
%\end{equation}
%The regularizating results also imply that, if the path  $\zeta$ is Brownian,  there are time intervals of fixed length migrating to infinity along which the stochastic viscosity solutions are actually $C^{1,1}$-in space for a set of times of Hausdorff dimension $1/2$.
%\smallskip
%
%

When $\zeta$ is in $C_0^1([0,\infty);\R)$ or $\text{BV}_0([0,\infty);\R)$,   in \eqref{takis1}  ``$\cdot$'' stands for multiplication and the problem  falls within the scope of the classical Crandall-Lions theory of viscosity solutions. When $\zeta \notin \text{BV}([0,\infty);\R)$,  \eqref{takis100} is studied using the notion of pathwise (stochastic) viscosity solutions, which was introduced and is been developed by the authors in  \cite{LS1, LS2, LS4, S}. In this setting ``$\cdot$'' only signifies the way the path, which can be nowhere differentiable, acts on $H$. When $\zeta$ is a Brownian motion, then ``$\cdot$'' is the usual ``$\circ$'' in the Stratonovich calculus. Note that  the solutions do not have sufficient regularity to actually interpret  the equation  in this sense. 
%\smallskip
Pathwise solutions of \eqref{takis100} are well posed in $\text{BUC}(\R^d \times [0,\infty))$,  %$\text{BUC}(\R^d \times [0,\infty))$ stands for 
the set of bounded uniformly continuous functions on $\R^d \times [0,\infty)$. 
\smallskip

The two ``deterministic'' regularity results are stated in terms of the symmetric matrix
$$F(p):=\sqrt{D^2H(p)}.$$
The first claim is about the regularizing effect of \eqref{takis1}.  We remark that all the inequalities and solutions below should be understood in the viscosity sense.
\begin{theorem}\label{takis21} Assume \eqref{hconvex2}. If $u\in \text{BUC}(\R^d \times [0,\infty))$
is a solution of $u_t=H(Du) $  
(resp. \\ $u_t=-H(Du)$) \ $\text{in $\R^d \times (0,\infty)$}$ and, for some $C \in (0,\infty]$,  
\begin{equation}\label{manos1}
-F(Du(\cdot, 0))D^2u(\cdot, 0) F(Du(\cdot, 0)) \leq CI \ \text{in} \ \R^d, 
\end{equation}
(resp.
 \begin{equation}\label{manos2}
-F(Du(\cdot, 0))D^2u(\cdot, 0) F(Du(\cdot, 0)) \geq - CI  \ \text{in} \ \R^d),
\end{equation} 
then, for all $t >0$, 
\begin{equation}\label{takis22}
 -F(Du(\cdot, t))D^2u(\cdot, t) F(Du(\cdot, t)) \leq \dfrac{C}{1+Ct} I  \ \text{in} \ \R^d,
 \end{equation}
 (resp.
\begin{equation}\label{takis23}
 -F(Du(\cdot, t))D^2u(\cdot, t) F(Du(\cdot, t)) \geq -\dfrac{C}{1+Ct} I  \ \text{in} \ \R^d).
 \end{equation}
 \end{theorem}
Estimates  \eqref{takis22} and \eqref{takis23}  are sharper versions of the classical regularizing effect-type estimates for viscosity solutions (see Lions \cite{lbook}, Lasry and Lions~\cite{ll}), which say that, if $u_t=H(Du)$  
(resp. $u_t=-H(Du)$) in $\R^d \times [0,\infty)$, and, for some $C \in (0,\infty]$, $-D^2u(\cdot, 0) \leq C I$ (resp. $-D^2u(\cdot, 0) \geq -C I) \ \text{in} \ \R^d$, then, for all $t >0$, 
\begin{equation}\label{takis24}
 -D^2 u(\cdot, t) \leq \dfrac{C}{1+\theta Ct} I  \ \text{in} \ \R^d
 \end{equation}
 (resp.
\begin{equation}\label{takis25}
 -D^2 u(\cdot, t)  \geq -\dfrac{C}{1+\theta Ct} I  \ \text{in} \ \R^d.)
 \end{equation}
 Note that, when $C=\infty$, that is,  no assumption is made on $u(\cdot,0)$, then \eqref{manos1} and \eqref{manos2} reduce to 
\beq\label{manos22.1}  -F(Du(\cdot, t))D^2u(\cdot, t) F(Du(\cdot, t)) \leq \dfrac{1}{t} \quad (\text{resp.} -F(Du(\cdot, t))D^2u(\cdot, t) F(Du(\cdot, t)) \geq -\dfrac{1}{t}),\eeq
%(resp. 
%\beq\label{takis23.1}-F(Du(\cdot, t))D^2u(\cdot, t) F(Du(\cdot, t)) \geq -\dfrac{1}{t}), \eeq
which are sharper versions of \eqref{takis24} and \eqref{takis25}, in the sense that they do not depend on $\theta$, of the classical 
estimates 
$$-D^2 u(\cdot, t) \leq \dfrac{1}{\theta t} \quad (\text{resp.} \  D^2 u(\cdot, t) \geq -\dfrac{1}{\theta t}).$$

\smallskip

We continue with  the propagation of regularity result by first recalling what was known. Indeed, it was shown in \cite{ll} that, if $u$ solves $u_t=H(Du)$ (resp. $u_t=-H(Du)$) with $H$ satisfying \eqref{hconvex1}, then,   
\begin{equation}\label{takis2221}\text{ if \ $-D^2 u(\cdot, 0) \geq - C I$, \ then \  $-D^2u(\cdot,t) \geq - \dfrac{C}{(1-\Theta Ct)_+}$},\end{equation}
(resp.
\beq\label{takis2222} \text{if  \ $-D^2u(\cdot, 0) \leq C I$, \ then \ $-D^2u(\cdot,t) \leq  \dfrac{C}{(1-\Theta Ct)_+}$.})\eeq

The new propagation of regularity result depends on the dimension. In what follows, we say that $H:\R^d\to \R$ is quadratic, if there exists a symmetric matrix $A$ which satisfies \eqref{hconvex1}  such that 
$$H(p)=(Ap,p).$$

\begin{theorem}\label{takis26.1} Assume \eqref{hconvex2} and let   $u\in \text{BUC}(\R^d \times [0,\infty))$ solve $u_t=H(Du) $   
(resp.   $u_t=-H(Du)$) \ $\text{in $\R^d \times (0,\infty)$}$.  Suppose that   $d=1$ or $H$ is quadratic.  If, for some  $C>0$, 
\begin{equation}\label{manos4}
 -F(Du(\cdot, 0))D^2u(\cdot, 0) F(Du(\cdot, 0)) \geq -CI \ \text{in} \ \R^d, 
 \end{equation}
 (resp. 
 \begin{equation}\label{manos5}
 -F(Du(\cdot, 0))D^2u(\cdot, 0) F(Du(\cdot, 0)) \leq CI  \ \text{in} \ \R^d), 
 \end{equation}
  then, for all $t >0$, 
\begin{equation}\label{takis27}
 -F(Du(\cdot, t))D^2u(\cdot, t) F(Du(\cdot, t)) \geq - \dfrac{C}{(1-Ct)_+} I  \ \text{in} \ \R^d,
 \end{equation}
 (resp.
\begin{equation}\label{takis28}
 -F(Du(\cdot, t))D^2u(\cdot, t) F(Du(\cdot, t)) \leq \dfrac{C}{(1-Ct)_+} I  \ \text{in} \ \R^d.)
 \end{equation}
\end{theorem}

The result for $d\geq 2$ and general $H$ is more restrictive.

\begin{theorem}\label{takis26} Let $d>1$ and assume that $H$ is not quadratic and satisfies  \eqref{hconvex2}. Let   $u\in \text{BUC}(\R^d \times [0,\infty))$ solve $u_t=H(Du) $   
(resp.  $u_t=-H(Du)$) \ $\text{in $\R^d \times (0,\infty)$}$ and assume that $u(\cdot,0) \in C^{1,1}(\R^d)$. If, for some  $C>0$, 
\begin{equation}\label{manos4}
 -F(Du(\cdot, 0))D^2u(\cdot, 0) F(Du(\cdot, 0)) \geq -CI \ \text{in} \ \R^d, 
 \end{equation}
 (resp. 
 \begin{equation}\label{manos5}
 -F(Du(\cdot, 0))D^2u(\cdot, 0) F(Du(\cdot, 0)) \leq CI  \ \text{in} \ \R^d), 
 \end{equation}
  then, for all $t >0$, 
\begin{equation}\label{takis27}
 -F(Du(\cdot, t))D^2u(\cdot, t) F(Du(\cdot, t)) \geq - \dfrac{C}{(1-Ct)_+} I  \ \text{in} \ \R^d,
 \end{equation}
 (resp.
\begin{equation}\label{takis28}
 -F(Du(\cdot, t))D^2u(\cdot, t) F(Du(\cdot, t)) \leq \dfrac{C}{(1-Ct)_+} I  \ \text{in} \ \R^d.)
 \end{equation}
\end{theorem}
\smallskip

It turns out that  the assumption that $u(\cdot,0) \in C^{1,1}(\R^d)$ if $d>1$ and $H$ is not quadratic   is necessary to have estimates like \eqref{takis27} and \eqref{takis28}. This is the claim  of  the next theorem.
\begin{theorem}\label{takis29} Assume \eqref{hconvex2} and  $d>1$.  If  \eqref{takis27} holds for all solutions $u \in \text{BUC}(\R^d \times [0,\infty))$ of $u_t=H(Du) $   
(resp.  $u_t=-H(Du)$) \ $\text{in $\R^d \times (0,\infty)$}$ with $u \in C^{0,1}(\R^d)$ satisfying 
\eqref{manos4} (resp. \eqref{manos5}), then the map  $\lambda \to (D^2H(p +\lambda \xi) \xi^\perp, \xi^\perp)$ must be concave (resp. convex).  In particular, both estimates hold  without any restrictions on the data if and only if 
$H$ is quadratic.
\end{theorem}

The motivation behind Theorem~\ref{takis21} and Theorem~\ref{takis26.1} and Theorem~\ref{takis26} is twofold. Firstly, we wish to obtain as sharp as possible regularity results for solutions of \eqref{takis1}. Secondly, we want to see if it is possible 
to obtain intermittent regularity results for  \eqref{takis100}, like the ones obtained in \cite{gg} in the specific case that $H(p)=\frac{1}{2}|p|^2$, where, of course, $\theta=\Theta=1$, 
$F(Du)D^2uF(Du)=D^2u$  and the ``new'' estimates are the same as the old ones, that is, \eqref{takis2221} and \eqref{takis2222},  which hold without any regularity conditions.  
\smallskip

The regularity results of \cite{gg} follow from an iteration of \eqref{takis24}, \eqref{takis25}, \eqref{takis2221} and \eqref{takis2222}. As we describe next,  the iteration scheme cannot, however,  work when $H$ is not quadratic unless $d=1$. To explain problem, we consider the first two steps of the possible iteration consider $u\in \text{BUC}(\R^d \times (0,\infty))$ which solves 
$$u_t=H(Du) \ \text{in} \ R^d\times (0,a], \quad  u_t=-H(Du) \ \text{in} \ R^d\times (a,a+b] \quad  \text{and} \quad u_t=H(Du) \ \text{in} \ R^d\times (a+b,a+b+c].$$
If the only estimates available were \eqref{takis24},  \eqref{takis25},  \eqref{takis2221} and  \eqref{takis2222}, we find, after a simple algebra, that 
$$D^2u(\cdot, a) \geq -\dfrac{1}{\theta a} I, \quad D^2 u(\cdot, a+b) \geq - \dfrac{1}{(\theta a - \Theta b)_+} I  \quad \text{and} \quad  D^2 u(\cdot, a+b+c) \geq - \dfrac{1}{(\theta a - \Theta b)_+ + \theta c} I.$$
It is immediate that the above estimates cannot be iterated unless there is a special relationship  between the time intervals and the convexity constants which will, something which not be possible for arbitrary continuous paths $\zeta$. 
\smallskip

We discuss next what would happen,  if it were possible to use the estimates of Theorem~\ref{takis26} without any regularity restrictions, as it is the case when $d=1$. To simplify the notation, we introduce the matrix 
\beq\label{ellie2}
{\mathcal W} (t):= F(Du(\cdot,t)) D^2u(\cdot,t) F(Du(\cdot,t)),
\eeq
and observe that, using  Theorem~\ref{takis21}, Theorem~\ref{takis26.1} and Theorem~\ref{takis26},   we would have
$${\mathcal W} (a) \geq -\dfrac{1}{a} I, \quad {\mathcal W} (a+b) \geq -\dfrac{1}{(a-b)_+} I \quad \text{and} \quad {\mathcal W} (a+b+c) \geq -\dfrac{1}{(a-b)_+ +c} I,$$
which can be further iterated, since the estimates are expressed only in terms of  increments $\zeta$.
\smallskip

Before we move to the intermittent regularity results, it is necessary to make some additional remarks. For the shake of definiteness, we continue the discussion in the context of the example above. Although $u(\cdot, a)$ may not be in $C^{1,1}$, it follows from \eqref{takis2221} and \eqref{manos5} that, for some $h\in (0,b]$ and $t\in (a, a+h)$, $u(\cdot, t)\in C^{1,1}$. There is no way, however, to guarantee that $h=b$. Moreover, as we show in section~2, in general,  it is possible to have $u$ and $h>0$ such that $u_t=-H(Du)$ in $\R^d\times (-h,0]$, $u_t=H(Du)$  in $\R^d\times (0,h]$, $u(\cdot,t) \in C^{1,1}$ for $t\in (-h,0) \cup (0,h)$ and $u(\cdot,0) \notin C^{1,1}$. The implication is that when $d>1$ and $H$ is not a quadratic, there is no hope to obtain after iteration smooth solutions. 
\smallskip

We continue with the discussion of the intermittent regularity.  To state the  results,  
 we introduce  the running maximum and minimum functions $M:[0,\infty)\to \R$ and $m:[0,\infty)\to \R$ of a path $\zeta \in 
C_0([0,\infty);\R)$ defined respectively by 
\begin{equation}\label{running}
M(t):=\underset{0\leq s \leq  t}\max \zeta(t) \quad \text{and} \quad m(t):=\underset{0\leq s \leq  t}\max \zeta(t).
\end{equation}
\begin{theorem}\label{takis310}  Assume   \eqref{hconvex2} and   $d=1$ or that $H$ is quadratic  when  $d>1$, fix  $\zeta \in C_0([0,T);\R)$ and let $u \in \text{BUC}(\R^d \times [0,\infty)$ be a solution of  \eqref{takis100}. 
Then, for all $t>0$, 
\begin{equation}\label{takis311}
- \dfrac{1}{M(t)-\zeta(t)}\leq -F(Du(\cdot,t))D^2u(\cdot,t)F(Du(\cdot,t)) \leq \dfrac{1}{\zeta(t)-m(t)}.  
\end{equation}
\end{theorem}
\smallskip

Note that when \eqref{takis311} holds, then, at times $t$ such that $m(t) < \zeta (t) < M(t)$, $u(\cdot, t)\in C^{1,1}(\R^d)$ with bounds 
independent of the Hamiltonian, since 
\eqref{takis311} implies that, for all $t>0$, 
\begin{equation}\label{IRR}
|F(Du(\cdot, t)) D^2u (\cdot, t))F(Du(\cdot, t))| \leq  \max\left[\frac{1}{\zeta(t) -m(t)}, \frac{1}{M(t) -\zeta(t)}\right].
\end{equation}
\smallskip

When, however, \eqref{takis311} is not available, the best regularity estimate available, which is also new, is a decay on the  Lipshitz constant $\|Du\|$ of;
in what follows $\|\cdot\|$ stands for the usual $L^\infty$-norm.
\begin{theorem}\label{takis32} Assume \eqref{hconvex2},  fix  $\zeta \in C_0([0,T);\R)$, let $u \in \text{BUC}(\R^d \times [0,\infty)$ be a solution of  \eqref{takis100}. Then, for all $t>0$,
\begin{equation}\label{takis33}
\|Du(\cdot,t)\| \leq \sqrt{\dfrac{2\|u(\cdot,t)\|}{\theta (M(t)-m(t))}}.
\end{equation}
\end{theorem}
\smallskip

It follows from \eqref{takis33} that, for any $t>0$ such that $m(t)<M(t)$, any solution of \eqref{takis100} is actually Lipshitz continuous.

\smallskip

An immediate consequence of the estimates in Theorem~\ref{takis32} and Theorem~\ref{takis310}, which is based on well known properties of the Brownian motion (see, for, example, Peres~\cite{p})  is the following observation.
\begin{theorem}\label{takis1000}
Assume that $\zeta$ is a Brownian motion and $H$ satisfies \eqref{hconvex2}. There exists an uncountable subset of $(0,\infty)$ with no isolated points and of Hausdorff measure $1/2$, which depends on $\zeta$, off of which,  any stochastic viscosity solution of \eqref{takis100} is in  satisfies $C^{0,1}(\R^d)$ with a bound satisfying \eqref{takis33}. If  $d=1$ or $H$ is quadratic, for  the same set of times,  the solution is in  $C^{1,1}(\R^d)$ and satisfies \eqref{takis311}. 
\end{theorem}
The long time behavior of the solutions of \eqref{takis100} is an immediate consequence of Theorem~\ref{takis32}.
\begin{theorem}\label{takis34}
Assume \eqref{takis712},  fix $\zeta \in C_0([0,T);\R)$, and let $u \in \text{BUC}(\R^d \times [0,\infty)$ be a space periodic solution of  \eqref{takis100}.  If there exists  $t_n\to \infty$ such that  $M(t_n)-m(t_n) \to \infty$, then there exists  $u_\infty \in \R$ such that, as $t\to\infty$ and uniformly in space, $u(\cdot,t) \to u_\infty$.
\end{theorem}
(resp. \eqref{takis28} ) holds. 
It turns out that the assumption on the path in the previous theorem is again a well known property of the Brownian paths; again see \cite{p}. Therefore we have the following result.

\begin{theorem}\label{takis9}
Assume \eqref{takis712}. For almost every Brownian path $\zeta$, if $u \in \text{BUC}(\R^d \times [0,\infty))$ is a periodic solution of \eqref{takis100}, 
there exists a constant $u_\infty =u_\infty (\zeta, u(\cdot, 0))$ such that, 
%if $u$ is the stochastic viscosity solution of \eqref{takis1}, then, 
as $t\to \infty$ and uniformly in $\R^d$, $u(\cdot, t) \to u_\infty$.  Moreover, the random variable is, in general, not constant. 
\end{theorem}

\subsection*{Organization of the paper} The paper is organized as follows. In the next section we prove Theorem~\ref{takis21}, Theorem~\ref{takis26} and Theorem~\ref{takis29}. In section 2 we review several facts from the theory of pathwise viscosity solutions. We prove the intermittent regularity results   in section 3 and  we discuss  the long time behavior in section 4.  

\section {The regularity results for \eqref{takis1}}

The aim here is to prove 
Theorem~\ref{takis21}, Theorem~\ref{takis26.1}, Theorem~\ref{takis26} and Theorem~\ref{takis29}. Since the proofs are similar, we only present them for 
\begin{equation}\label{takis400}
u_t=H(Du) \ \text{in} \ \R^d\times [0,\infty).
\end{equation}
The convexity of $H$ allows for the use of  the classical  
Lax-Oleinik formula which yields (see \cite{lbook}) that, for any $t > 0$ and $L$ the convex conjugate of $H$,  
\begin{equation}\label{takis20}
u(x,t)=\underset{ y \in \R^d} \sup [u(y, 0) - t \ L(\dfrac{x-y}{t})].
\end{equation}
It is immediate  that  \eqref{takis20} can be rewritten as 
\begin{equation}\label{takis2000}
u(x,t)= \sup [u(y,0) - tL(\dfrac{z}{t}): \ y, z \in \R^d \ \text{and} \ y+z=x]. 
\end{equation}

Since $u_0$ is bounded and $H$ is uniformly convex, it follows that, for each $(x,t) \in \R^d\times (0,\infty)$, the 
maximum in \eqref{takis2000} is achieved at pairs  $(\bar y, \bar z)=(\bar y(x,t), \bar z (x,t)),$  
that is, there exist $\bar y, \bar z \in \R^d$ such that 
\begin{equation}\label{takis2001}
u(x,t)=u_0(\bar y) - tL(\dfrac{\bar z}{t}) \quad \text{and} \quad \bar y + \bar z= x.
\end{equation}

For future use, we remark that, given $\chi \in \R^d$ and $\lambda \in [0,1]$, it follows from \eqref{takis2000} and \eqref{takis2001} that 
\begin{equation}\label{takis2002}
u(x\pm \chi ,t) \geq u(\bar y(x,t),0) \pm (1-\lambda) \chi) - t L(\dfrac{\bar z(x,t) \pm \lambda \chi}{t}).
\end{equation}

As mentioned in the introduction,   the Lax-Oleinik formula has a regularizing effect and propagates regularity. The former property is that, for any $u (\cdot, 0)\in \text{BUC} (\R^d)$ and $t>0$, $u(\cdot,t)$ becomes semiconvex with an estimate that depends on $\theta$ in  \eqref{hconvex1} . The latter property has two parts.  The first is that the semiconvexity  $u(\cdot, 0)$ is propagated. The second is more delicate.  It says that, if $u(\cdot, 0)$ is semiconcave, 
then so does the solution up to a time that depends on the bound at $t=0$ and $\Theta$ in \eqref{hconvex1}.  Recall  that $u$ is semiconvex (resp. semiconcave) with constant $C$ if $D^2u \geq -CI \ (\text{ resp.} \  D^2u \leq CI).$
 
\smallskip

We summarize these facts in the next lemma. For a proof we refer to \cite{ll}. 

\smallskip
\begin{lemma}\label{takis2411} Assume \eqref{hconvex2} and let $u\in \text{BUC}(\R^d)\times [0,\infty)$ be a solution of \eqref{takis400}. Then:

(i)~If $u(0,\cdot)$ is semiconvex with constant $C$, 

 then, for all $t>0$,  
\begin{equation}\label{takis241}
D^2u(\cdot,t)  \geq -\dfrac{C}{1+\theta C} {I}. % \quad (\text{resp.} \ D^2S^+(t)u_0 \leq C{I}).
\end{equation}
(ii)~If $u(\cdot,0)$ is semiconcave with constant $C$, then  
\begin{equation}\label{takis242}
D^2 u(\cdot,t) \leq \dfrac{C}{(1-C\Theta t)_+}{I}.
\end{equation}
\end{lemma}
As discussed earlier the goal is to improve the bounds above in the sense that they become independent of the convexity constants of $H$. This requires to identify to study the propagation of the matrix $\mathcal W$, which was defined in \eqref{ellie2}. 
\smallskip

The reason is that, as it follows easily from the elementary calculations below, if the solution of \eqref{takis400}  is smooth, then ${\mathcal{ W}}$ satisfies the matrix equation
\begin{equation}\label{takis771}
{\mathcal W}_t= (DH(Du)) D{\mathcal{ W}} + {\mathcal{ W}}^2.
\end{equation} 
Indeed, when $u$ is smooth, it is easy to check that, since we are dealing with symmetric matrices,  $F(Du)$ and $D^2H$ satisfy the matrix equations 
$$F(Du)_t=DH(Du) D[F(Du)] \quad \text{and} \quad (D^2u)_t=DH(Du) D(D^2u) + D^2H(Du) (D^2u)^2.$$
It follows that 
$${\mathcal W}_t= DH D{\mathcal W} + FD^2H (D^2u)^2F,$$
and, since, in view of the symmetry,  $FD^2H (D^2u)^2F=FD^2uD^2HD^2uF=FD^2uFFD^2uF={\mathcal W}^2$, the claim follows.
\smallskip

It is then immediate, at least formally,  that ${\mathcal{ \overline W}}(t):=\max_x {\mathcal{ W}}(x,t)$ and ${\mathcal{ \underline W}}(t):=\min_x {\mathcal{W}}(x,t)$ satisfy respectively the differential inequalities
\begin{equation}\label{takis772}
\dot {\mathcal {\overline W}} \leq {\mathcal{ \overline W}}^2 \quad \text{and} \quad    \dot {\mathcal {\underline W}}
\geq {\mathcal{ \underline W}}^2,
\end{equation}
which easily yield, when everything is smooth, \eqref{takis22} and \eqref{takis27}.  

\begin{proof}[The proof of Theorem~\ref{takis21}]
We only show \eqref{takis22}. Although the conclusion follows when $u$ is smooth by the argument above, the claim does not rely on any regularity. Hence, it is necessary to come up with a different argument, which is based on the Lax-Oleinik formula. 
\smallskip

To simplify the notation,  we write $v$ and $v_0$ for  $u(\cdot, t)$ and $u(\cdot, 0)$ respectively, in which case  we have 
$$v(x)= \sup [v_0(y) -t L(\dfrac{z}{t}): y+z=x].$$ 
The claim is that, for any smooth function $\phi$,  any maximum $x_0$ of $v-\phi$ and all $\xi \in \R^d$,
$$-\left(F(D\phi(x_0))D^2\phi(x_0)F(D\phi(x_0))\xi, \xi \right)\leq \dfrac{C}{1+Ct} |\xi|^2.$$
Let  $y_0, z_0 \in \R^d$ be such that 
$y_0+z_0=x_0 \quad \text{and} \quad v(x_0)=v(y_0)-t L(\dfrac{z_0}{t}).$
\smallskip

It follows from Lemma~\ref{takis2411}, the assumption on $u(\cdot,0)$ and \eqref{hconvex2} that  both $v$ and $v_0$ are semiconvex.  Then a standard argument in the theory of viscosity solutions yields that without loss of generality we may assume that, up to a small linear perturbation,  $v$ and $v_0$ are respectively differentiable at $x_0$ and $y_0$ and 
$Dv(x_0)=Dv_0(y_0)=DL(\dfrac{z_0}{t})=D\phi(x_0)$, the latter coming from the Lax-Oleinik formula.   
\smallskip

Next fix $\xi \in \R^d$, $\lambda \in [0,1]$, set $\chi=F(D\phi(x_0))\xi$, and recall  \eqref{takis2002}. 

It follows that, for all $h>0$, 
$$ v_0(y_0 \pm h(1-\lambda)\chi) -t L(\dfrac{z_0 \pm h \lambda \chi}{t}) -\phi (y_0 \pm h(1-\lambda)\chi + z_0 \pm h \lambda \chi) \leq v(y_0)-t L(\dfrac{z_0}{t}) - \phi (x_0).$$
The choice of $\chi$ and the fact that $D^2L(q)=(D^2H(DL(q)))^{-1}$ yields  
$$-(F(D\phi(x_0))D^2\phi)(x_0)F(D\phi(x_0))\xi,\xi) \leq  [(1-\lambda)^2  +  \dfrac{\lambda^2}{t}]|\xi|^2.$$
Optimizing over $\lambda$ yields  the claim.

\end{proof}

\begin{proof}[The proof of Theorem~\ref{takis26}] We first  assume that  $u(\cdot,0) \in C^2(\R^d)$ in which case the solution of \eqref{takis400} can be constructed by the methods of characteristics up  to some time $T^\star$, which depends only on bounds on $D^2u(\cdot,0)$ and $D^2H$.
\smallskip

The conclusion follows from the equation satisfied by $\mathcal W$. Alternatively, we may use   the characteristics to find that 
$$Du(x,t)=Du(X^{-1}(x,t)) \ \ \text{and} \ \ 
D^2u(x,t)=D^2u_0(X^{-1}(x,t)) \dfrac{\partial }{\partial x} X^{-1}(x,t)$$
and 
$$\quad \dfrac{\partial }{\partial x} X^{-1}(x,t)=(I -t DH^2(Du_0(X^{-1}(x,t))D^2u_0(X^{-1}(x,t))^{-1}.$$
A simple computation now leads to \eqref{takis27} as long as the characteristics are invertible. 
It follows from the  convexity of the Hamiltonian and the estimates \eqref{takis22} and  \eqref{takis27} that the characteristics do not intersect up to $T^\star=1/C$.
\smallskip

To prove the general result  we approximate $u(\cdot, 0)\in C^{1,1}(\R^d)$ using a standard mollifier  $\rho_\ep$, that is, we take $u_0^\ep=u(\cdot,0)\star \rho_\ep,$ and consider the  initial value problem
$$u^\ep_t=H(Du^\ep) \ \ \text{in} \ \ \R^d\times (0,\infty) \quad u^\ep(\cdot, 0)=u_0^\ep  \ \ \text{in} \ \ \R^d.$$
It is clear that, as $\ep \to 0$,
 $u^\ep \to u$ uniformly on $\R^d\times [0,T]$ for every $T>0$. The claim follows in the limit $\ep \to 0$ by the stability property of the viscosity inequality, if we show a suitable upper bound for $\ep>0$.
\smallskip

Notice that $u^\ep$ satisfies \eqref{takis27} with a constant $C_\ep \to C$ as $\ep\to 0$. Moreover, the convolution operation preserves the $C^{1,1}$-bounds of $u_0$. Hence the interval of invertibility of the characteristics remains the same.

\end{proof}

We continue with the proof of the propagation of regularity result when $d=1$.
\begin{proof}[The proof of Theorem~\ref{takis26.1}]  
Let $G(p):=(H''(p))^{-1} $ and observe that the claim  we are interested in proving can be written as  
$$\text{if} \ -u_{xx}(\cdot, 0) \geq - CG(u_x(\cdot,0)), \ \text{then for all $t>0$,}  \ -u_{xx}(\cdot, t) \geq -\dfrac{G(u_x(\cdot, t) )C}{(1-Ct)_+}.$$

In view of the stability of the viscosity solutions, we assume that $H$ is smooth, and to circumvent the issues of regularity at $t=0$, we use the viscosity approximation 
\begin{equation}\label{takis511}
u^\ep_t -\ep H''(u^\ep_x) u^\ep_{xx} -H(u^\ep_x)=0 \ \text{in} \ \R\times (0,\infty) \quad u^\ep(\cdot,0)=u_0
\end{equation}
where $u_0$ satisfies  $-u_{0, xx}(\cdot, 0) \geq - CG(u_{0,x}(\cdot,0))$. 
\smallskip

For the rest of the proof we omit the dependence of $u^\ep$ and, for $h:[0,\infty) \to \R$ to be chosen,  we consider the auxiliary function $w(\cdot,t):=u_{xx} - h(t) G(u_x)$.
\smallskip

A straightforward computation yields 
\begin{equation*}
\begin{split}
w_t-\ep H''(u_x) w_{xx}-H'(u_x)w_x=H''(u_x) u^2_{xx} - h'G(u_x) +\ep H''''(u_x) u^3_{xx}\\
 +\ep h (2H'''(u_x) G'(u_x) + H''(u_x)G''(u_x)) u^2_{xx}.
\end{split}
\end{equation*}
The choice of $G$ simplifies the equation above to read
\begin{equation}\label{takis6000}
w_t -\ep H''(u_x)w_{xx} -H'(u_x)w_x= H''(u_x)(u^2_{xx} -h' G^2(u_x)) + \ep H''''(u_x)u^2_{xx} w.
\end{equation}
Select $m$ to be the solution of the ode  $h'=h^2$ with $m(0)=C$, that is $h(t)=\dfrac{C}{(1-Ct)_+}.$
\smallskip

Then evaluating \eqref{takis6000} the first positive time that $w$ achieves a maximum yields $w_t \leq 0$.
It follows that $w$ achieves its maximum at $t=0$, which, in view of the choice $m$ leads to
$w(\cdot, t) \leq 0$, and, hence, the result.  

\end{proof}

When $d\geq 2$,  the above proof amounts to considering,  for some appropriately chosen  uniformly elliptic symmetric matrix $A=A(p)$,  the viscous approximation
$u^\ep_t-\ep \text{tr} [A(Du^\ep) D^2 u^\ep]=H(Du^\ep)$
and to proving that 
$$\dfrac{d}{dt} \underset{x\in \R^d, \ |\xi|=1} \sup \left(F(Du^\ep) D^2u^\ep F(Du^\ep) \xi, \xi \right)
\leq \big [ \underset{x\in \R^d, \ |\xi|=1} \sup \left (F(Du^\ep) D^2u^\ep F(Du^\ep) \xi,\xi \right) \big]^2.
$$
Unfortunately the usual maximum principle-type arguments requires to assume signs for expressions involving $D^4H$ and $D^4L$, and, hence, do not provide any useful information.
\smallskip

The difficulty described above is, however, to expected in view of the Theorem~\ref{takis29} whose proof we present next. 

\begin{proof}[The proof of Theorem~\ref{takis29}] 
To avoid unnecessary complications we take  $d=2$ and consider a Hamiltonian   $H:\R^2 \to \R$  which satisfies  \eqref{hconvex2}.  
\smallskip

Without loss of generality we may assume that $u_0$ is not $C^1$ at the origin, and, after further reductions, we take
\beq\label{takis5002}u_0(x,y)=|x|-\frac{1}{2}y^2,\eeq
and 
consider  the initial value problem 
\begin{equation}\label{takis5001}
u_t + H(Du)=0 \ \text{in} \  \R^2\times (0,\infty) \ \text{and} \  u(\cdot,0)=u_0.
\end{equation}
After  translating and  rotating, we may reduce to that case that, for some $c_0>0$, 
\begin{equation}\label{takis5000}
 H(0,0)=0, DH(0,0)=0, \ \text{and} \ D^2H(0,0)=c_0 I. %  \ \text{and}  \ D^2H(\pm 1, 0) =I.
\end{equation}
To further simplify the presentation, we assume  that $H$ is even, that is, 
\begin{equation}\label{Heven}
H(-p)=H(p) \ \text{for all} \ p\in\R^2,
\end{equation}
and, in addition, that 
$$D^2H(\pm,0)=I.$$ 
The general case, which will leave it up to the reader, is to assume that $u_0$ is $C^1$ in the direction of $\xi\in \R^2$ and not $C^1$ in the direction $\xi^\perp$ and to consider  more general  Hamiltonians. 
\smallskip

Next we observe that, due to the finite speed of propagation  of the Hamilton-Jacobi equations, it is enough to work in a neighborhood of the origin and show that it is not possible 
to have in this case a propagation of regularity property like \eqref{takis28}.
\smallskip

In view of  \eqref{takis5000}, \eqref{Heven},  and the choice of $u_0$, in a small neighborhood of the origin and in the viscosity sense, we have
\begin{equation}\label{takis5003}
-F(Du_0)D^2u_0F(Du_0)= - I \ \begin{bmatrix}
    0       & 0\\ 
    0       & -1 \\
\end{bmatrix} \ I \leq I.
\end{equation}
Assume next that it is possible to have the lower bound claimed in \eqref{takis28}. Then,   in a neighborhood of the origin,  in the viscosity sense and for $t>0$, we must have 
\begin{equation}\label{takis800.1}
-F(Du(\cdot,t))D^2u(\cdot,t)F(Du(\cdot,t))\leq \dfrac{1}{(1-t)_+} I.
\end{equation} 
The regularizing effect of \eqref{takis5001}, it  follows that,  for  $t\in (0,1-\ep)$, $u(\cdot,t)\in C^{1,1}(\R^2)$. 
Moreover, since $u_0$ is even, we conclude that $Du(0,0,t)=0$ in $(0,1).$ Finally, the $C^{1,1}$ regularity of $u$ in space, the last observation and the facts that $H(0,0)=0$ and $u_0(0,0)=0$ imply that $u(0,0,t)=0$ in $(0,1).$
\smallskip

Next we perform a second-order blow up of $u$ at $0$, that is, we consider
\[u_\ep(x,t):=\dfrac{1}{\ep}u(\sqrt {\ep} x, \sqrt {\ep}y,t),\]
which solves 
\begin{equation}\label{takis9141}
u_{\ep, t} =  H_\ep(Du_\ep) \ \ \text{in}  \ \ \R^2\times (0,1) \quad \quad u_\ep(\cdot,0)=u_{\ep, 0},
\end{equation}
with
$$H_\ep (q):=\dfrac{1}{\ep}  H(\sqrt{\ep} q) \quad \text{and} \quad u_{0, \ep}(x):=\dfrac{1}{\ep}u_0(\sqrt {\ep} x).$$
In view of the properties of $H$ at $p=0$,  it follows that 
$$H_\ep(q)=\int_0^1\int_0^1 (D^2H(\sqrt{\ep}\lambda q)q, q) \sigma d\lambda d\sigma,$$
and, hence, as $\ep \to 0$ and locally uniformly in $q$,
$$H_\ep (q) \to \frac{1}{2} (D^2H(0,0) q, q)=\frac{1}{2} c_0 |q|^2.$$
Set ${\mathds{1}}_{\{0\}} := 0 \ \text{if \  $x=0$ \ and $\infty \ \text{otherwise},$}$. Then, as $\ep \to 0$ and locally uniformly in $\R\setminus \{0\} \times \R$, 
$$u_{0,\ep}(x,y) \to {\mathds{1}}_{\{0\}} -\frac{1}{2}y^2.$$
Since the scaling preserves the second-derivative bound, $u(0,0,t)=0$ and $Du(0,0,t)=0$,   for $t\in (0,1)$, the $u^\ep$'s are uniformly bounded in a neighborhood of $(0,0)$. Finally, as $\ep\to 0$,  and, hence,  have  second derivatives which are bounded uniformly in $\ep$. 
\smallskip

It  follows that, in a neighborhood of the origin, the $u^\ep$'s converge to $v$ which solves 
$v_t +  \frac{1}{2} c_0 |Dv|^2=0 \ \text{in} \ \R^2\times (0,1)$.  Using the Lax-Oleinik formula for u and, hence,$u_\ep$, we find that 
$$v(x,y,t)=\frac{1}{2tc_0}x^2 -\frac{1}{2(1-c_0 t)}y^2,$$
 which is, in view of the results of Crandall, Lions and Souganidis~\cite{cls},
 the unique solution of 
\begin{equation}\label{takis900}
v_t +  \frac{1}{2} c_0 |Dv|^2=0 \ \text{in} \ \R^2\times (0,1) \quad v(x,y,0)= {\mathds 1}_{\{0\}}(x) -\frac{1}{2}y^2.
\end{equation}
Next we remark that the desired estimate is stable under limits. Hence,
if, at (x,y)=(0,0), we had $\sqrt{D^2H(0,0)}D^2u(0,0)\sqrt{D^2H(0,0)} \leq -1/(1-t) I$ as assumed,
then, in the blow up limit $\ep\to 0$, we must have
\[ -c_0 \begin{bmatrix}
    \dfrac{1}{c_0 t}       & 0\\ 
    0       & -\dfrac{1}{1-c_0t} \\
\end{bmatrix} \leq \frac{1}{1-t} I,\] 
which requires that $c_0 \leq 1.$ 
\smallskip

Note that $c_0\leq 1$ also implies that $$D^2H( \dfrac{1}{2} (-1,0) +  \dfrac{1}{2} (1,0))_{22}=[D^2H(0,0)]_{22}=c_0 \leq I =\dfrac{1}{2} [D^2H((-1,0))]_{22} +  \dfrac{1}{2}  [D^2H((1,0))]_{22},$$
which yields the convexity assumption asserted in the claim for $p=0$ and $\xi=(1,0)$.
\smallskip

Finally, for the propagation of regularity result to be true for both signs in \eqref{takis1} for any $u(\cdot,0)$, it is necessary to have that map $t\to (D^2H(p+t\xi)\xi^\perp\xi^\perp)$ is at most linear for all $p,\xi$. Since, however, it is assumed that $D^2H$ is bounded, it must be the case that $D^2H$ is constant, which, in turn, implies that $H$ is a homogeneous quadratic.

\end{proof}

We conclude with an observation that adds to the statement of the lack of propagation of regularity without additional assumptions on $u$. The example in the proof above yields an initial datum for which there cannot be propagation of regularity. It may, however, be argued that such a $u$ cannot arise in the the process of the iteration. 
\smallskip

We show next that this is not the case. Indeed consider the solution $u$ of \eqref{takis5001} with $u_0$ given by \eqref{takis5002}. As it was agued above, $u \in C^{1,1} (\R^d\times (0,T))$ for some fixed $T>0$. Let $v(x,t)=u(x,-t)$. It is immediate that  $v \in C^{1,1} (\R^d\times (-T,0))$ and
$v_t - H(Dv)=0 \ \text{in} \ \R^d\times (-T,0).$

\section{Pathwise viscosity solutions}

We provide a brief overview of the theory of pathwise viscosity solutions. Instead of stating the intrinsic definition of the solution, here we work with the fact, which was  proven in the references mentioned earlier,   that the solution operator is obtained as the unique extension of the one defined for smooth paths using the classical Crandall-Lions theory of viscosity solutions. For details, we refer the readers to 
\cite{LS1, LS2, LS4, S}. 
\smallskip

We recall that the ``classical''  theory viscosity solution theory applies to initial value problems
\begin{equation}\label{takis10}
u_t=H(Du)\dot \zeta  \ \ \text{in} \ \  \R^d\times (0,\infty) \qquad u(\cdot, 0)=u_0 \ \ \text{in } \ \ \R^d,
\end{equation}

with $H\in C(\R^d), \ u_0 \in \text{BUC}(\R^d)$ and  $\zeta \in \text{BV}_0([0,T])$ for all $T>0$. 
\smallskip

When necessary to emphasize the dependence on the path, we write $u^\zeta$ for the solution of \eqref{takis10}.

\smallskip

The basic result of the ``deterministic'' theory is stated next.

\begin{theorem}\label{takis11}
Assume that $H\in C(\R^d)$. For each $T>0$, $\zeta \in \text{BV}_0([0,T])$ and $u_0 \in \text{BUC}(\R^d),$
the initial value problem \eqref{takis10} has a unique solution $u\in \text{BUC}(\R^d \times [0,\infty)).$ Moreover, if 
$u, v \in \text{BUC}(\R^d \times [0,\infty))$ are respectively sub- and super-solutions of \eqref{takis10}, then, for all $t\in [0,T]$, 
$$\|(u(\cdot, t)-v(\cdot, t))_+\| \leq \|(u(\cdot, 0)-v(\cdot, 0))_+\|.$$
Finally, if $u_0\in C^{0,1}(\R^d)$, then,  for all $t\in [0,T]$, $u(\cdot, t) \in C^{0,1}(\R^d)$ and $\|Du(\cdot,t)\|\leq \|Du(\cdot, 0)\|.$
\end{theorem}

The main result about the pathwise solutions of \eqref{takis1} is stated next; for the proof we refer to~\cite{LS1, LS2, LS4, S}.

\begin{theorem}\label{takis12}
For any $\zeta \in C_0([0,\infty);\R)$ and $u_0 \in \text{BUC}(\R^d),$  the initial value problem \eqref{takis1} has a unique pathwise (stochastic) solution if and only if $H$ is the difference of two convex functions. In addition, the contraction and Lipschitz continuity properties in Theorem~\ref{takis11} are also true.
\end{theorem} 

We continue with a summary of the key properties of the pathwise solutions that we will be using in the paper;
for proofs we refer again to~\cite{S}.

\begin{proposition}\label{takis13}
(i)~The pathwise solutions of \eqref{takis1} are continuous with respect to the Hamiltonian $H$, the path $\zeta$,  and the initial value $u_0$.

(ii)~Assume that $H$ is the difference of two convex functions and consider families $(\zeta_\ep)_{\ep\in (0,1)}, (\zeta_{\ep'})_{\ep'\in (0,1)} \in  C_0^1([0,\infty))$ and $(u_{0,\ep})_{\ep\in (0,1)}, (u_{0,\ep'})_{\ep'\in (0,1)} \in \text{BUC}(\R^d)$  such that, as $\ep, \ep' \to 0$, $\zeta_\ep -  \zeta_{\ep'} \to 0$ locally uniformly in $[0,\infty)$ and $u_{0,\ep}-u_{0,\ep'} \to 0$ uniformly in $\R^d$. Let $u_\ep, u_{\ep}' \in \text{BUC}(\R^d \times [0,\infty))$ be  the viscosity solutions of 
\eqref{takis10} with paths $\zeta_\ep$, and $\zeta_{\ep'}$ and initial datum $u_{0,\ep}$, and $u_{0,\ep'}$ respectively. Then, as  $\ep, \ep' \to 0$ and uniformly in $ \text{BUC}(\R^d \times [0,\infty))$, for all $T>0$, $u_{\ep} - u_{\ep'} \to 0$. Moreover, if, as $\ep, \ep' \to 0, \ \zeta_\ep, \zeta_{\ep'} \to \zeta$ and $u_\ep, u_{\ep '} \to u$, then $u$ is the pathwise solution of \eqref{takis1} with path $\zeta$ and initial datum $u(\cdot, 0).$
\end{proposition}

The next two facts are  a consequence of the contraction property.

\begin{proposition}\label{takis14}
Assume that $H$ is the difference of two convex functions, $H(0)=0$,  and let  $u \in \text{BUC}(\R^d \times [0,\infty))$ be a  pathwise solution of \eqref{takis100}.

(i)~The maps $t \to  \sup_{\R^d}  u(\cdot,t)$ and $t \to  \inf_{\R^d} u(\cdot,t)$  
are respectively nonincreasing and nondecreasing.

(ii)~If $u(\cdot,0) \in C^{0,1}(\R^d)$, then, for all  $t\in  (0,\infty))$, $u(\cdot,t) \in C^{0,1}(\R^d)$ and the map $t \to \|Du(\cdot, t)\|$ is nonincreasing.

\end{proposition}
The last item discussed here  is the recent of work of Gassiat, Gess and the authors~\cite{ggls}, who used the so-called skeleton function of a path, to obtain finite speed of propagation and domain of dependence-type results for the pathwise solutions of \eqref{takis10}, The representation of the solutions that comes from the skeleton plays an important role in the proof of the intermittent regularity.
\smallskip

Given $\zeta \in C_0([0,T])$, the sequence $(\tau_{i})_{i\in\mathbb{Z}}$ of successive extrema of $\xi$  is defined by 
\begin{equation}
\tau_{0}:=\sup\left\{ t\in[0,T]: \ \zeta(t)=M(t)\mbox{ or }\zeta(t)=m(t)\right\} ,\label{eq:extrema1}
\end{equation}
and, for all  $ i\geq0$, \begin{equation}
\tau_{i+1}=\left\{ \begin{array}{ll}
\arg\max_{[\tau_{i},T]}\zeta & \mbox{ if }\ \ \zeta(\tau_{i})<0,\\[1.5mm]
\arg\min_{[\tau_{i},T]}\zeta & \mbox{ if } \ \ \zeta(\tau_{i})>0.
\end{array}\right.\label{eq:extrema2}
\end{equation}
The skeleton  path  ${R}_{0,T}(\zeta)$  of $\zeta \in C_0([0,T])$ is 
 the  piecewise linear function agreeing  with $\zeta$ on $(\tau_{i})_{ i\in\mathbb{Z}}$.
\smallskip

The usefulness of the skeleton is seen in the next result, which was shown in \cite{ggls}. In what follows, we write $u^\zeta$ to denote the solutions of \eqref{takis400}. 
\begin{theorem}\label{takis1000}
Fix $\zeta \in C_0([0,T])$ and consider its  skeleton path ${R}_{0,T}(\zeta)$. Then,
\begin{equation}\label{takis1001}
u^\zeta(\cdot, T)=u^{{R}_{0,T}(\zeta)} (\cdot, T).
\end{equation}
\end{theorem}

\section{Intermittent regularity}

To discuss the two regularizing results, that is, Theorem~\ref{takis310} and Theorem~\ref{takis32},  it is necessary to introduce some additional notation, which also explains the method of the proof.

\smallskip

In what follows we denote by by $S^\pm$ the solution operators  of the initial value problems 
\begin{equation}\label{takis2010}
u^\pm_t=\pm H(Du^\pm) \ \  \text{in} \ \ \R^d\times[0,\infty) \qquad u^\pm(\cdot,0)=u_0 \ \ \text{in} \ \ \R^d,
\end{equation}
that is, for $u_0\in \text{BUC}(\R^d)$ and $t\in (0,\infty)$, $S^\pm(t)u_0 \in  \text{BUC}(\R^d)$ is the solution of \eqref{takis2010}.
\smallskip

\begin{proof}[The proof of Theorem~\ref{takis32}] Fix $\zeta\in C_0([0,\infty)$, a solution $u^\zeta \in  \text{BUC}(\R^d) \times [0,\infty)$ of \eqref{takis10} and $t>0$. In view of its definition, the up to $t$ skeleton path ${R}_{0,t}(\zeta)$ of $\zeta$ contains an interval of length $M(t)-m(t)$.  Since, in view of the construction of the skeleton, $\tau_0$ is a time at which $M(t)$ is achieved, $m(t)$ is achieved either at $\tau_1$, in which case we have
$$u^{{R}_{0,t}(\zeta)} (\cdot, \tau_1)=S^-(M(t)-m(t))u^{{R}_{0,t}(\zeta)} (\cdot, \tau_0),$$
or at $\tau_{-1}$ and, hence,
$$u^{{R}_{0,t}(\zeta)} (\cdot, \tau_0)=S^+(M(t)-m(t))u^{{R}_{0,t}(\zeta)} (\cdot,\tau_{-1}).$$
It follows that, in the first  case,
$$-D^2u^{{R}_{0,t}(\zeta)} (\cdot, t_1)\geq - \dfrac{1}{\theta(M(t)-m(t))} I,$$
while in the second case
$$-D^2u^{{R}_{0,t}(\zeta)} (\cdot, t_1)\leq \dfrac{1}{\theta (M(t)-m(t))} I.$$
Next we note that, in view of decreasing  in time  property of  the $\|\cdot\|_\infty$-norm of the solutions and \eqref{takis1001}, we  clearly have $\|u^{{R}_{0,t}(\zeta)}\|_\infty \leq \|u\|_\infty$.
\smallskip

A standard estimate yields in the first case (resp. second) case that 
$$\|Du^{{R}_{0,t}(\zeta)}(\cdot,\tau_1)\|\leq \sqrt{\dfrac{2\|u\|_\infty}{\theta (M(t)-m(t))}}, \quad 
(resp. \ 
\|Du^{{R}_{0,t}(\zeta)}(\cdot,\tau_0)\|\leq \sqrt{\dfrac{2\|u\|_\infty}{\theta (M(t)-m(t))}}.)$$
Since the Lipshitz constant decreases in time either of the two estimates above imply that
$$\|Du^{{R}_{0,t}(\zeta)}(\cdot,t)\|\leq \sqrt{\dfrac{2\|u\|_\infty}{\theta (M(t)-m(t))}},$$
and the claim now follows in view of \eqref{takis1001}.
 
 \end{proof}

\smallskip

Next, fix $\zeta \in C^1_0([0,\infty))$ and $u_0\in \text{BUC}(\R^d)$,  let $u\in \text{BUC}(\R^d \times [0,\infty))$ be the solution of \eqref{takis10} and assume 
that there exists a sequence $(t_n)_{n\in \N\cup \{0\}}$ in $[0,\infty)$  such that 
\begin{equation}\label{takisg1}
\begin{cases}
t_0=0 \ \ \text{and} \  \ \dot \zeta (t_n)=0 \  \text{ for $n\in \N$,  \ and }\\[2mm] 
\text{$\dot \zeta>0$  \ in \ $(t_{2k},t_{2k+1})$ \ and \ $\dot \zeta <0$ \ in \ $(t_{2k+1},t_{2k+2})$} \ \  \text{for $k\in \N\cup \{0\}.$}
\end{cases}
\end{equation}

 Note that we assume that $\dot \zeta =0$ only along the sequence $(t_n)_{n\in \N}$, since, if $\dot \zeta =0$ in some interval, $u$ is constant in time there. 
\smallskip

It is then immediate that
\begin{equation}\label{eth1}
u(\cdot, t)=S^+(\zeta(t)-\zeta(t_{2k}))u(\cdot, t_{2k}) \quad \text{if} \quad  t\in (t_{2k}, t_{2k+1}),
\end{equation}
and 
\begin{equation}\label{eth2}
u(\cdot, t)=S^-(\zeta(t_{2k+1})-\zeta(t)))u(\cdot, t_{2k+1}) \quad  \text{if} \quad t \in (t_{2k+1}, t_{2(k+1)}).
\end{equation}
\vskip.075in
We use now the above in order to prove the intermittent regularity result when $d=1$, that is, 
Theorem~\ref{takis310}. 
\begin{proof}[The proof of Theorem~\ref{takis310}] 
We show only the upper bound. Since the claim for a general continuous path follows by density,  we assume that  $\zeta$ is smooth.
\smallskip

The argument is based on iterating  \eqref{takis27} and \eqref{takis28} and using  the substitutions in \eqref{eth1} and \eqref{eth2}.
\smallskip

We make next the above precise. We fix a partition as in \eqref{takisg1}.  It is then immediate from \eqref{takis27} that 
\begin{equation}\label{takisx}
F(u_x(\cdot, t_1)u_{xx}(\cdot, t_1) F(u_x(\cdot,t_1)) \geq - \dfrac{1}{\zeta(t_1)-\zeta(0)}\ I\geq - \dfrac{1}{M(t_1)-\zeta(t_1)} \ I.
\end{equation}
Having established \eqref{takisx},  we now proceed with the general argument using induction. For simplicity we assume  that the lengths of the partition intervals are such that the bounds never blow up. 
\smallskip

Assume next that we know the result up to $t_n$ and prove it in $(t_n,t_{n+1}).$ If $n=2k$,  the induction
claim gives 
Then, for $t \in (t_{2k}, t_{2k+1})$, $u(\cdot, t)=S^+(\zeta(t)-\zeta(t_{2k}))u(\cdot, t_{2k})$.  It follows 
from Theorem~\ref{takis26.1} that
\[ \begin{split} F(u_x(\cdot, t))u_{xx}(\cdot, t) F(u_x(\cdot, t) )& \geq - \dfrac{ \dfrac{1}{M(t_{2k})-\zeta(t_{2k})}}{{1+\dfrac{1}{M(t_{2k})-\zeta(t_{2k})}(\zeta(t)-\zeta(t_{2k}))}}I\\[3mm]
& \geq - \dfrac{1}{M(t_{2k})-\zeta(t_{2k}) + \zeta(t)-\zeta(t_{2k})} I 
 \geq - \dfrac{1}{M(t) -\zeta(t)} I.
\end{split}
\]
\smallskip

If $t \in (t_{2k+1}, t_{2(k+1)})$, then $u(\cdot, t)=S^-(\cdot, \zeta(t_{2k+1})-\zeta(t))u(\cdot, t_{2k+1})$, then repeating the argument that led us to the first step of the iteration we find
$$F(u_x(\cdot, t))u_{xx}(\cdot, t) F(u_x(\cdot, t))  \geq - \dfrac{\dfrac{1}{M(t_{2k+1}) -\zeta(t_{2k+1})}}{\left(1-\dfrac{\zeta(t_{2k+1})-\zeta(t_{2(k+1)})}{M(t_{2k+1}) -\zeta(t_{2k+1})}\right)_+} I \geq - \dfrac{1}{M(t) -\zeta(t)} I.
$$
\end{proof}
\smallskip

Theorem~\ref{takis1000} follows from Theorom~\ref{takis32} and Theorem~\ref{takis310} 
and the properties of the Brownian motion-see, for example, \cite{p}.
\smallskip

\section{The long time behavior}

We consider solutions of \eqref{takis1} and are interested in their behavior as $t\to \infty$. In order to avoid technicalities due to the behavior of the solutions $u\in \text{BUC}(\R^d)$ at infinity,  throughout this section we work with periodic functions in ${\mathbb {T}}^d$.

To explain the problem we first look at two very simple cases. In the first case, we fix some $V\in \R^d$ and consider first the linear initial value
problem
$$du=(V, Du) \cdot d\zeta \ \ \text{in} \ \  \R^d \times (0,\infty) \quad u(\cdot,0)=u_0 \ \ \text{in} \ \ \R^d.$$

Its solution is
$u(x,t)=u_0( x + V\zeta(t)),$
and clearly we cannot expect that $u(\cdot,t)$ has, as $t\to \infty$, a uniform limit.
\smallskip

We also look at \eqref{takis1} with $H$ satisfying \eqref{takis712} and $\dot\zeta >0$ and $\underset{t \to \infty}\lim \zeta(t)=\infty$.
Since\\  
$u(x,t)=\underset{y\in \R^d} \sup [u_0(y) -tL(\dfrac{x-y}{\zeta(t)})],$
it is immediate that, as $t\to \infty$ and uniformly in $x$, $u(x,t) \to \sup u$.

We prove next the main asymptotic result of the paper when $H$ satisfies \eqref{takis712}. 

\smallskip

\begin{proof}[The proof of Theorem~\ref{takis34}]   The contraction property and the fact that $H(0)=0$ yield that  the family  $(u(\cdot,t))_{t >0}$ is uniformly bounded. 
\smallskip

We assume next that the Hamiltonian satisfies \eqref{hconvex2}. It then follows from the intermittent regularizing  property \eqref{takis33}, the assumption  on $M-m$ and the fact that  the Lipschitz constant of the solutions decreases in time that, as $t\to \infty$,  $\|Du(\cdot,t)\| \to 0$. 
\smallskip

In view of the periodicity, we find that, along subsequences $s_n\to \infty$, the $u(\cdot,s_n)$ converge uniformly to constant. 

\smallskip

It remains to show that the whole family converges to the same constant. This is consequence of \eqref{takis33}, the periodicity and Proposition~\ref{takis14}, which yields that 

\begin{equation}\label{takisxx}
t \to \max_{x \in \R^d} u(x,t) \quad \text{is non-increasing, and}  \quad t \to \min_{x \in \R^d} u(x,t) \quad \text{is non-decreasing}.
\end{equation}
It remains to remove the assumption that the Hamiltonians satisfy \eqref{hconvex2}.  Indeed, If   \eqref{takis712} holds , we approximate $H$ uniformly by a sequence $(H_m)_{m\in \N}$ of Hamiltonians satisfying \eqref{hconvex2}. Let $u_m$ be the solution of the \eqref{takis1} with Hamiltonian $H_m$ and same initial datum. Since, as $m\to infty$,  $u_m \to u$ uniformly in $\R^d \times [0,\infty)$ for all $T>0$, it follows that, for all $t>0$, 
$$\int_{\mathbb{T}} H(Du(x,t))  dx \leq \underset{m \to \infty} \liminf \int_{\mathbb{T}} H(D u_m (x,t)) dx.$$

We may now use the sequence $t_n$ as before to conclude.

\qed

\smallskip

We conclude with an example that shows that, in the stochastic setting,  the limit constant $u_\infty$ may be a random variable. 
\smallskip

We consider the initial value problem
\begin{equation}\label{1}
du=|u_x| \circ dB \ \ \text{in} \ \ \R\times (0,\infty) \quad u(\cdot,0)=u_0 \ \ \text{in} \ \ \R^d,
\end{equation}
where $u_0$ is periodic  (with period $2$) on $\R$ and on $[0,2]$,  $u_0(x)=1-|x-1|.$ 

Let $c$ be the limit as $t \to \infty$ of $u$. Since $1-u_0(x)=u(x+1)$ and $-B$ is also a Brownian motion with the same law as $B$, it follows that 
\begin{equation}\label{2}
\mathcal L(c)=\mathcal L(1-c),
\end{equation}
where $\mathcal L(f)$ denotes the law of the random variable $f$.
\smallskip

If the limit $c$ of the solution of \eqref{1} is deterministic, then  \eqref{2} implies that $c\equiv 1/2$.  We show next that this is not the case.

\smallskip

Recall that the pathwise solutions are Lipshitz with the respect to paths in the sense that if $u_1,u_2$ are two pathwise solutions of \eqref{1} with paths respectively $\zeta_1, \zeta_2$, then there exists $L>0$, which depends on $\|u_{0,x}\|$  such that, for any $T>0$,
\begin{equation}\label{3}
\underset{x\in \R, t \in  [0,T]}\max |u_1(x,t)-u_2(x,t)| \leq L \underset{t \in  [0,T]} \max|\zeta_1(t)-\zeta_2(t)|.
\end{equation}

Next we fix $T=2$ and  use \eqref{3}  to compare the solutions of \eqref{1} with $\zeta_1\equiv B$ and $\zeta_2(t)=t$ and  $\zeta_1\equiv B$ and $\zeta_2(t)=-t$.
\smallskip

When $\zeta_2=t$ (resp. $\zeta=-t$) the solution $u_2$  of \eqref{1} is given by 
$$u(x,t)=\underset{|y|\leq t}\max u_0(x+y) \quad (\text{resp.} \ \  
u(x,t)=\underset{|y|\leq t}\min u_0(x+y)).$$   

It  is then simple to check that, if $\zeta_2(t)=t$, then $u_2(\cdot, 2)\equiv 1,$ while, when $\zeta_2(t)=-t$, $u_2(x,2)=0$. 
\smallskip

Fix $\ep=1/4L$ and consider the events
\begin{equation}\label{11}
A_+:=\{\underset{t\in [0,2]} \max|B(t) -t| < \ep\} \quad \text{and} \quad A_-:=\{\underset{t\in [0,2]} \max|B(t) + t| < \ep\}.
\end{equation}
Of course, 
\begin{equation}\label{12}
\mathbb{P}(A_+) >0 \quad \text{and} \quad \mathbb{P}(A_-) >0.
\end{equation}

Then  \eqref{3} implies 
\begin{equation}\label{takis12}
u_1(x,2) \geq 1-L\ep = 3/4  \ \ \text{on} \ \ A_+ \quad \text{and} \quad u_2(x,2) \leq L\ep = 1/4   \ \ \text{on} \ \ A_-.
\end{equation}

It follows that the random variable $c$  cannot be constant since in $A_+$ it must be bigger than $3/4$ and in $A_-$ smaller than $1/4$.

\smallskip

In an upcoming publication (Gassiat, Lions and Souganidis~\cite{gls}) we are visiting this problem and obtain in a special case more information about $u_\infty$.

\end{proof}

\bigskip

\noindent ($^{1}$) Coll\`{e}ge de France,
11 Place Marcelin Berthelot, 75005 Paris, 
and  
CEREMADE, 
Universit\'e de Paris-Dauphine,
Place du Mar\'echal de Lattre de Tassigny,
75016 Paris, FRANCE\\ 
email: lions@ceremade.dauphine.fr
\\ \\
\noindent ($^{2}$) Department of Mathematics 
University of Chicago, 
5734 S. University Ave.,
Chicago, IL 60637, USA, \ \  
email: souganidis@math.uchicago.edu
\\ \\
($^{3}$) Partially supported by the Air Force Office for Scientific Research  grant FA9550-18-1-0494.
\\ \\ 
($^{4}$)  Partially supported by the National Science Foundation grant DMS-1600129, the Office for Naval Research grant N000141712095 and the Air Force Office for Scientific Research grant FA9550-18-1-0494.

\end{document}